\documentclass[11pt]{article}

\usepackage{amsmath,amssymb,theorem}
\usepackage[mathscr]{eucal}
\usepackage[all]{xy}

\newtheorem{theorem}{Theorem}[section]

\newtheorem{proposition}[theorem]{Proposition}

{\theorembodyfont{\rmfamily}
\newtheorem{example}[theorem]{Example}}

\makeatletter \@addtoreset{equation}{section} \makeatother

\title{\bf QUANTUM GROUPS AND BOUNDED SYMMETRIC DOMAINS}

\author{\tt S. Sinel'shchikov \and\tt  L. Vaksman}
\date{Mathematics Department, Institute for Low Temperature Physics and
Engineering,\\ 47 Lenin Ave, Kharkov 61103, Ukraine}

\begin{document}
\large

\maketitle

\begin{abstract}
Recent results of the authors on quantum bounded symmetric domains and
quantum Harish-Chandra modules are expounded.
\end{abstract}

\stepcounter{section}

{\bf\thesection.} Here is a general outline of problems in non-compact
quantum group theory to be discussed in this talk.

The geometrical methods of studying Harish-Chandra modules are very
important in representation theory of real reductive groups. Our obstacle
that arises in passage from ordinary to quantum groups is that necessary
methods of non-commutative geometry are not at a hand. This forces us to
combine the process of searching for geometric realizations of quantum
Harish-Chandra modules with looking for suitable concepts and results of
non-commutative geometry.

For such approach to be successful, one has to choose properly the class of
modules to be considered. Our choice is based on the notion of quantum
bounded symmetric domain \cite{SV}.

It is well known that every irreducible bounded symmetric domain admits a
standard realization in a complex prehomogeneous vector space. We begin
with a result on $*$-algebras of `polynomials on quantum prehomogeneous
vector spaces'.

Turn to precise formulations.

\bigskip

\stepcounter{section}

{\bf\thesection.} We assume $\mathbb{C}$ as a ground field. The parameter
$q$ is supposed to be in the interval $(0,1)$. All the algebras are unital,
unless the contrary is stated explicitly.

Consider the $*$-algebra $\operatorname{Pol}(\mathbb{C})_q$ given by a
single generator $z$ and the defining relation
\begin{equation*}\label{z_comm}
z^*z=q^2zz^*+1-q^2.
\end{equation*}
$\operatorname{Pol}(\mathbb{C})_q$ is a quantum analogue of the polynomial
algebra on a complex plane treated as a two dimensional real vector space.
We are interested in representations of \hbox{$*$-algebras} by bounded
linear operators.

Let $T_F$ be a representation of $\operatorname{Pol}(\mathbb{C})_q$ in
$l^2(\mathbb{Z}_+)$ given by
\begin{equation}
T_F(z)e_n=\left(1-q^{2(n+1)}\right)^\frac12e_{n+1},\quad T_F(z^*)e_n=
\begin{cases}
\left(1-q^{2n}\right)^\frac12e_{n-1}, & n\in\mathbb{N},
\\ 0, & n=0,
\end{cases}
\end{equation}
with $\{e_n\}$ being the standard basis of $l^2(\mathbb{Z}_+)$. This is the
so called Fock representation.

One can readily use the spectral theory of operators and the commutation
relations
\begin{equation}
zy=q^{-2}yz,\qquad z^*y=q^2yz^*,
\end{equation}
with $y=1-zz^*$, to derive the following easy and well known

\begin{proposition}\label{C}
\begin{enumerate}
\item $T_F$ is a faithful irreducible $*$-representation of
$\operatorname{Pol}(\mathbb{C})_q$ by bounded linear operators in a Hilbert
space.

\item Every representation with these properties is unitarily equivalent to
$T_F$.
\end{enumerate}
\end{proposition}

This Proposition could be treated as a q-analogue of the Stone-von Neumann
theorem since the commutation relation between the elements
$$a^+=(1-q^2)^{-\frac12}z,\qquad a=(1-q^2)^{-\frac12}z^*$$
is a q-analogue of the canonical commutation relation $aa^+-a^+a=1$.

During 90's q-analogues were found for polynomial algebras on some special
prehomogeneous vector spaces. This became a background for producing an
advanced extension of Proposition \ref{C} and lead to the notion of quantum
bounded symmetric domain. Describe a construction of those q-analogues for
polynomial algebras.

Consider a simple complex Lie algebra $\mathfrak{g}$ of rank $l$ whose
Cartan matrix is $a=(a_{ij})$. Up to an isomorphism this Lie algebra can be
described in terms of the generators $\{H_i,E_i,F_i\}_{i=1,2,\ldots,l}$ and
the well known defining relations. The linear span $\mathfrak{h}$ of
$H_1,H_2,\ldots,H_l$ is a Cartan subalgebra, and the linear functionals
$\alpha_1,\alpha_2,\ldots,\alpha_l$ on $\mathfrak{h}$ defined by
$$\alpha_j(H_i)=a_{ij},\qquad i,j=1,2,\ldots,l,$$
form a system of simple roots for the Lie algebra $\mathfrak{g}$.

Choose a simple root $\alpha_{l_0}$ which appears in the decomposition
$\delta=\sum\limits_{i=1}^ln_i\alpha_i$ of the maximal root $\delta$ with
coefficient 1.\footnote{Such simple roots exist for all simple complex Lie
algebras except $E_8$, $F_4$, $G_2$.} Let $H\in\mathfrak{h}$ be the element
given by
$$\alpha_{l_0}(H)=2,\qquad \alpha_j(H)=0,\quad j\ne l_0,$$
and
$$
\mathfrak{k}=\left\{\xi\in\mathfrak{g}|\:[H,\xi]=0\right\},\qquad
\mathfrak{p}^\pm=\left\{\xi\in\mathfrak{g}|\:[H,\xi]=\pm 2\xi\right\}.
$$
Then one has
$\mathfrak{g}=\mathfrak{p}^-\oplus\mathfrak{k}\oplus\mathfrak{p}^+$, which
constitutes a $\mathbb{Z}$-gradation:
$$
[\mathfrak{k},\mathfrak{k}]\subset\mathfrak{k},\qquad
[\mathfrak{k},\mathfrak{p}^\pm]\subset\mathfrak{p}^\pm,\qquad
[\mathfrak{p}^+,\mathfrak{p}^-]\subset\mathfrak{k},\qquad
[\mathfrak{p}^+,\mathfrak{p}^+]=[\mathfrak{p}^-,\mathfrak{p}^-]=0.
$$
Consider the complex simply connected Lie group with the Lie algebra
$\mathfrak{g}$ and its subgroup $K$ of those elements which preserve the
gradation
$$
K=\{g\in
G|\:\operatorname{Ad}_g\mathfrak{k}=\mathfrak{k},\;
\operatorname{Ad}_g\mathfrak{p}^\pm=\mathfrak{p}^\pm\}.
$$
We follow H. Rubenthaler in calling $\mathfrak{p}^\pm$ the prehomogeneous
vector spaces of commutative parabolic type. The prehomogeneity means
existence of an open $K$-orbit.

A construction of the $*$-algebra $\operatorname{Pol}(\mathfrak{p}^-)_q$, a
quantum analogue of the polynomial algebra on $\mathfrak{p}^-$, can be
found in \cite{SV}. Describe the outline of this construction and one of
its steps.

For that, we need a Hopf $*$-algebra $(U_q\mathfrak{g},*)$, with
$U_q\mathfrak{g}$ being the quantum universal enveloping algebra. This Hopf
$*$-algebra is given by its generators
$$K_i,K_i^{-1},E_i,F_i,\qquad i=1,2,\ldots,l,$$
the well known Drinfeld-Jimbo relations \cite{Jant}. The comultiplication
$\triangle$, the counit $\varepsilon$, the antipode $S$, and the involution
$*$ are given by
$$
\triangle(E_i)=E_i\otimes 1+K_i\otimes E_i,\qquad \triangle(F_i)=F_i\otimes
K_i^{-1}+1\otimes F_i,\qquad \triangle(K_i)=K_i\otimes K_i,
$$
$$\varepsilon(E_i)=\varepsilon(F_i)=0,\qquad\varepsilon(K_i)=1,$$
$$S(E_i)=-K_i^{-1}E_i,\qquad S(F_i)=-F_iK_i,\qquad S(K_i)=K_i^{-1},$$
$$
(K_j^{\pm 1})^*=K_j^{\pm 1},\;j=1,2,\ldots,l;
$$
$$
E_j^*=
\begin{cases}
K_jF_j, & j\ne l_0
\\ -K_jF_j, & j=l_0
\end{cases};\qquad
F_j^*=
\begin{cases}
E_jK_j^{-1}, & j\ne l_0
\\ -E_jK_j^{-1}, & j=l_0
\end{cases}.
$$
In the special case $l=1$ one has the Hopf $*$-algebra
$U_q\mathfrak{su}_{1,1}\overset{\mathrm{def}}{=}(U_q\mathfrak{sl}_2,*)$.

In what follows all the $U_q\mathfrak{g}$-modules are assumed to be weight:
$$
V=\bigoplus_{\lambda=
(\lambda_1,\lambda_2,\ldots,\lambda_l)\in\mathbb{Z}^l}V_\lambda,\qquad
V_\lambda=\{v\in V|\:K_i^{\pm 1}v=q_i^{\pm\lambda_i}v,\;i=1,2,\ldots,l\},
$$
with $q_i=q^{d_i}$ and $d_i$ being the coprime numbers that symmtrize the
Cartan matrix: $d_ia_{ij}=d_ja_{ji}$. Define the linear operators
$H_1,H_2,\ldots,H_l$ in $V$ by $H_i|_{V_{\boldsymbol{\lambda}}}=\lambda_i$.
Obviously, in $\operatorname{End}V$ one has the relations
$$K_i^{\pm 1}=q_i^{\pm H_i},\qquad i=1,2,\ldots,l.$$

Consider the Hopf subalgebra $U_q\mathfrak{k}$ generated by
$$K_j^{\pm 1},\quad j=1,2,\ldots,l;\qquad E_i,F_i,\quad i\ne l_0.$$
The weight finitely generated $U_q\mathfrak{g}$-module $V$ is called a
Harish-Chandra module if it splits as a sum of simple finite dimensional
$U_q\mathfrak{k}$-modules, each of those having a finite multiplicity in
$V$.

The construction and classification of simple quantum Harish-Chandra
modules are important open problems.

Turn back to a construction of the $*$-algebra
$\operatorname{Pol}(\mathfrak{p}^-)_q$. Its relation to quantum bounded
symmetric domains is to be based on the fact that
$\operatorname{Pol}(\mathfrak{p}^-)_q$ is a $(U_q\mathfrak{g},*)$-module
algebra. Recall the latter notion.

Consider an algebra $F$ which is also a module over a Hopf algebra $A$. $F$
is called an $A$-module algebra if the multiplication $m:F\otimes F\to F$
is a morphism of $A$-modules. In the case of a unital algebra $F$ one has
to require additionally an $A$-invariance of its unit. In the presence of
involutions in $A$ and in $F$ they have to be compatible:
$$(af)^*=(S(a))^*f^*,\qquad a\in A,\;f\in F.$$

The initial step towards $\operatorname{Pol}(\mathfrak{p}^-)_q$ is in
producing a $U_q\mathfrak{g}$-module algebra
$\mathbb{C}[\mathfrak{p}^-]_q$, which is a q-analogue for the algebra
$\mathbb{C}[\mathfrak{p}^-]$ of holomorphic polynomials on
$\mathfrak{p}^-$. We follow the idea of V. Drinfeld for producing by
duality function algebras on quantum groups. To construct the
$U_q\mathfrak{g}$-module algebra $\mathbb{C}[\mathfrak{p}^-]_q$, we start
with considering the `dual' coalgebra. This coalgebra is going to be a
generalized Verma module, specifically the $U_q\mathfrak{g}$-module with a
generator $v^{(0)}$ and defining relations
$$
E_iv^{(0)}=(K_i^{\pm 1}-1)v^{(0)}=0,\quad i=1,2,\ldots,l;\qquad
F_jv^{(0)}=0,\quad j\ne l_0.
$$
The comultiplication is essentially determined by
$$\triangle v^{(0)}=v^{(0)}\otimes v^{(0)}$$
(see \cite{SV} for details). It is worthwhile to note that the
multiplicities of simple finite dimensional $U_q\mathfrak{k}$-modules in
$\mathbb{C}[\mathfrak{p}^-]_q$ are the same as the multiplicities of
corresponding simple finite dimensional $U\mathfrak{k}$-modules in
$\mathbb{C}[\mathfrak{p}^-]$.

In the special case $l=1$ one has the polynomial algebra in one variable
$z$ and with the following $U_q\mathfrak{sl}_2$-action:
\begin{equation}\label{KF_act}
K^{\pm 1}f(z)=f(q^{\pm 2}z),\qquad Ff(z)=
q^\frac12\,\frac{f(q^{-2}z)-f(z)}{q^{-2}z-z},
\end{equation}
\begin{equation}\label{E_act}
Ef(z)=-q^\frac12z^2\,\frac{f(z)-f(q^2z)}{z-q^2z}.
\end{equation}

The next step is in producing the $U_q\mathfrak{g}$-module algebra
$\mathbb{C}[\mathfrak{p}^+]_q$ of `antiholomorphic polynomials on the
quantum vector space $\mathfrak{p}^-$'. It is a very easy step:
\begin{itemize}
\item[-] $\mathbb{C}[\mathfrak{p}^+]_q$ is just
$\mathbb{C}[\mathfrak{p}^-]_q$ as an Abelian group,

\item[-] the identity map
\begin{equation}\label{invol}
*:\mathbb{C}[\mathfrak{p}^-]_q\to\mathbb{C}[\mathfrak{p}^+]_q,\qquad
*:f\mapsto f
\end{equation}
is antilinear,

\item[-] the action of $U_q\mathfrak{g}$ on $\mathbb{C}[\mathfrak{p}^+]_q$
is given by
$$
(\xi f)^*=(S(\xi))^*f^*,\qquad\xi\in
U_q\mathfrak{g},\;f\in\mathbb{C}[\mathfrak{p}^-]_q.
$$
\end{itemize}
The algebras $\mathbb{C}[\mathfrak{p}^\pm]_q$ are equipped with gradations
as follows:
$$
\mathbb{C}[\mathfrak{p}^\pm]_q=
\bigoplus_{j=0}^\infty\mathbb{C}[\mathfrak{p}^\pm]_{q,\mp j},\qquad
\mathbb{C}[\mathfrak{p}^\pm]_{q,j}=\{f\in\mathbb{C}[\mathfrak{p}^\pm]_q|\:
Hf=2jf\}.
$$
One can demonstrate that those algebras are generated by their homogeneous
components $\mathbb{C}[\mathfrak{p}^\pm]_{q,\mp 1}$ and they are quadratic
algebras.

At the final step the vector space
$$
\operatorname{Pol}(\mathfrak{p}^-)_q\overset{\mathrm{def}}{=}
\mathbb{C}[\mathfrak{p}^-]_q\otimes\mathbb{C}[\mathfrak{p}^+]_q
$$
is equipped with a structure of $(U_q\mathfrak{g},*)$-module algebra. An
involution $*$ is defined via the antilinear map \eqref{invol} in an
obvious way. A multiplication is imposed via the Drinfeld universal
R-matrix. Specifically, introduce the notation
$$
m^\pm:\mathbb{C}[\mathfrak{p}^\pm]_q\otimes\mathbb{C}[\mathfrak{p}^\pm]_q\to
\mathbb{C}[\mathfrak{p}^\pm]_q,\qquad m^\pm:f_1\otimes f_2\mapsto f_1f_2
$$
for multiplications in $\mathbb{C}[\mathfrak{p}^\pm]_q$ and a q-analogue
for the permutation of tensor multiples
$$
\check{R}:\mathbb{C}[\mathfrak{p}^+]_q\otimes\mathbb{C}[\mathfrak{p}^-]_q\to
\mathbb{C}[\mathfrak{p}^-]_q\otimes\mathbb{C}[\mathfrak{p}^+]_q,
$$
determined by the action of the universal R-matrix in
$\mathbb{C}[\mathfrak{p}^+]_q\otimes\mathbb{C}[\mathfrak{p}^-]_q$.

$\check{R}$ is well defined, as one can demonstrate that the weights of the
$U_q\mathfrak{g}$-module $\mathbb{C}[\mathfrak{p}^-]_q$ are non-negative
while the weights of the $U_q\mathfrak{g}$-module
$\mathbb{C}[\mathfrak{p}^+]_q$ are non-positive (the zero weight determines
one-dimensional weight spaces $\mathbb{C}\cdot 1$).

A multiplication
$m:\left(\mathbb{C}[\mathfrak{p}^-]_q\otimes\mathbb{C}[\mathfrak{p}^+]_q
\right)^{\otimes
2}\to\mathbb{C}[\mathfrak{p}^-]_q\otimes\mathbb{C}[\mathfrak{p}^+]_q$ is
defined via $m^\pm$ and $\check{R}$:
$$
m\overset{\mathrm{def}}{=}(m^-\otimes
m^+)\left(\operatorname{id}_{\mathbb{C}[\mathfrak{p}^-]_q}\otimes\check{R}
\otimes\operatorname{id}_{\mathbb{C}[\mathfrak{p}^+]_q}\right).
$$

Such multiplication equips $\operatorname{Pol}(\mathfrak{p}^-)_q$ with a
structure of $(U_q\mathfrak{g},*)$-module algebra \cite{SV}.

\begin{example}
Let $l=1$. One has $Ez^*=q^{-\frac32}$,\ \ \ \ \ $q^{-\frac{H\otimes
H}2}z^*\otimes z=q^2z^*\otimes z$,
$$
R(z^*\otimes z)=(1+(q^{-1}-q)E\otimes F)q^{-\frac{H\otimes H}2}z^*\otimes
z=q^2z^*\otimes z+1-q^2,
$$
$$\check{R}(z^*\otimes z)=q^2z\otimes z^*+1-q^2,\qquad z^*z=q^2zz^*+1-q^2.$$
\end{example}

A description of $\operatorname{Pol}(\mathfrak{p}^-)_q$ in terms of
generators and relations in a more general case
$\mathfrak{g}=\mathfrak{sl}_{\,l+1}$ is given in \cite{SSV}. A complete
list of irreducible $*$-representations is known for some of those
$*$-algebras \cite{PW, Turowska}.

Note that $\mathbb{C}[\mathfrak{p}^+]_{q,-1}\hookrightarrow
\mathbb{C}[\mathfrak{p}^+]_q\hookrightarrow
\operatorname{Pol}(\mathfrak{p}^-)_q$. A vector $v_0\ne 0$ from a space of
$*$-representation $T$ of $\operatorname{Pol}(\mathfrak{p}^-)_q$ is called
a vacuum vector if $T(f)v_0=0$ for all
$f\in\mathbb{C}[\mathfrak{p}^+]_{q,-1}$.

\begin{theorem}
\begin{enumerate}
\item There exists a unique (up to a unitary equivalence) faithful
irreducible $*$-representation $T_F$ of
$\operatorname{Pol}(\mathfrak{p}^-)_q$ by bounded linear operators in a
Hilbert space $H_F$.

\item There exists a unique (up to a scalar multiple) vacuum vector for
$T_F$ as above.
\end{enumerate}
\end{theorem}

To match correspondence with the special case
$\mathfrak{g}=\mathfrak{sl}_2$ considered above, we keep the term `Fock
representation' for $T_F$.

\bigskip

\stepcounter{section}

{\bf\thesection.} Let $d\nu$ be an invariant measure on $\mathbb{D}$. Our
goal is to obtain a q-analogue for the algebra $\mathscr{D}(\mathbb{D})$ of
smooth functions on $\mathbb{D}$ with compact supports and a q-analogue for
the invariant integral
$$
\nu:\mathscr{D}(\mathbb{D})\to\mathbb{C},\qquad
\nu:f\mapsto\int\limits_\mathbb{D}fd\nu.
$$
These are to be used to introduce a q-analogue for the space $L^2(d\nu)$
which is involved into formulating the principal problem of harmonic
analysis in $\mathbb{D}$. (Note that polynomials are not in $L^2(d\nu)$:
$\int\limits_{\mathbb{D}}1d\nu=\infty$.)

Consider the space of the Fock representation $T_F$ and the one-dimensional
orthogonal projection $P_0$ onto the vacuum subspace. To produce
$\mathscr{D}(\mathbb{D})_q$, it would be well to have an element
$f_0\in\operatorname{Pol}(\mathfrak{p}^-)_q$ such that
\begin{equation}\label{f0_P0}
T_Ff_0=P_0,
\end{equation}
but it does not exist. This is a motivation to attach such $f_0$ to
$\operatorname{Pol}(\mathfrak{p}^-)_q$.

Let us act formally. Extend the $(U_q\mathfrak{g},*)$-module algebra
$\operatorname{Pol}(\mathfrak{p}^-)_q$ by attaching an element $f_0$ which
satisfies the following relations that are motivated by \eqref{f0_P0} and a
continuity argument:
\begin{itemize}
\item[-] $f_0^2=f_0$,\ \ $f_0^*=f_0$,

\item[-] $\psi^*f_0=f_0\psi=0$ for all weight
$\psi\in\mathbb{C}[\mathfrak{p}^-]_q$, $\psi\notin\mathbb{C}\cdot 1$,

\item[-] $K_j^{\pm 1}f_0=f_0$,\ \ $j=1,2,\ldots,l$,

$F_jf_0=
\begin{cases}
-\frac{q_{l_0}^{1/2}}{q_{l_0}^{-2}-1}f_0z_\mathrm{low}^*, & j=l_0
\\ 0, & j\ne l_0
\end{cases}$, \ \ \ \ \ $E_jf_0=
\begin{cases}
-\frac{q_{l_0}^{1/2}}{1-q_{l_0}^2}z_\mathrm{low}f_0, & j=l_0
\\ 0, & j\ne l_0.
\end{cases}$
\end{itemize}
Here $z_\mathrm{low}$ is the unique element of
$\mathbb{C}[\mathfrak{p}^-]_q$ with the properties
$$
K_i^{\pm 1}z_\mathrm{low}=q_i^{\pm a_{il_0}}z_\mathrm{low},\qquad
F_iz_\mathrm{low}=
\begin{cases}
q_{l_0}^\frac12, & i=l_0
\\ 0, & i\ne l_0
\end{cases}.
$$

We keep the notation $T_F$ for the natural extension of the Fock
representation onto the above algebra. The two-sided ideal
$\mathscr{D}(\mathbb{D})_q$ of this algebra generated by $f_0$ will be
called the algebra of finite functions on the quantum bounded symmetric
domain $\mathbb{D}$. $\mathscr{D}(\mathbb{D})_q$ is a
$(U_q\mathfrak{g},*)$-module algebra.

Recall that a linear functional $\nu$ on an $A$-module algebra $F$ is
called an invariant integral if $\nu(af)=\varepsilon(a)\nu(f)$ for all
$a\in A$, $f\in F$.

\begin{proposition}
There exists a non-zero $U_q\mathfrak{g}$-invariant integral on the
$U_q\mathfrak{g}$-module algebra $\mathscr{D}(\mathbb{D})_q$. It is unique
up to a constant multiple and can be chosen to be positive:
$\int\limits_{\mathbb{D}_q}f^*fd\nu>0$, with $f\ne 0$.
\end{proposition}

The integral looks like a q-trace. We formulate this result. Since
$\mathscr{D}(\mathbb{D})_q=\mathbb{C}[\mathfrak{p}^-]_qf_0
\mathbb{C}[\mathfrak{p}^+]_q$, one has
$\mathcal{H}_F\overset{\mathrm{def}}{=}\mathscr{D}(\mathbb{D})_qf_0=
\mathbb{C}[\mathfrak{p}^-]_qf_0=\operatorname{Pol}(\mathfrak{p}^-)_qf_0$.

Obviously, $\mathcal{H}_F$ is a $\mathscr{D}(\mathbb{D})_q$-module, a
$\operatorname{Pol}(\mathfrak{p}^-)_q$-module, and a weight
\hbox{$U_q\mathfrak{k}$-module}. Introduce the notation $\mathscr{T}_F$ for
the corresponding representations of $\mathscr{D}(\mathbb{D})_q$ and
$\operatorname{Pol}(\mathfrak{p}^-)_q$ in the vector space $\mathcal{H}_F$
(the Hilbert space $H_F$ is a completion of the pre-Hilbert space
$\mathcal{H}_F$).

\begin{proposition}
Let $\rho$ be the half sum of positive roots,
$\rho=\frac12\sum\limits_{i=1}^ln_i\alpha_i$ and
$\check{\rho}=\frac12\sum\limits_{i=1}^ln_id_iH_i$. The linear functional
on $\mathscr{D}(\mathbb{D})_q$
\begin{equation}\label{int_explicit}
\int\limits_{\mathbb{D}_q}fd\nu=\mathrm{const}\cdot
\operatorname{tr}\left(\mathscr{T}_F(f)q^{-2\check{\rho}}\right),\qquad
\mathrm{const}>0,
\end{equation}
is a positive $U_q\mathfrak{g}$-invariant integral.
\end{proposition}

Now \eqref{int_explicit} can be readily used to obtain a q-analogue of
weighted Bergman spaces and the well known geometric realizations for the
so called holomorphic discrete series in those spaces.

\bigskip

\stepcounter{section}

{\bf\thesection.} In the classical case $q=1$ one has a well known
geometric realization for the non-degenerate principal series of
Harish-Chandra modules on the open $K$-orbit $\Omega$ of the flag variety.

It is easy to obtain a quantum analogue for such a geometric realization,
and then to use it as a background for producing a quantum analogue for the
principal series. We restrict ourselves to producing a
$U_q\mathfrak{g}$-module algebra $\mathbb{C}[\Omega]_q$, which is a quantum
analogue for the algebra $\mathbb{C}[\Omega]$ of regular functions on the
affine algebraic variety $\Omega$.

Recall a general background on spherical weights. Every
$\lambda=(\lambda_1,\lambda_2,\ldots,\lambda_l)\in\mathbb{Z}_+^l$
determines a simple finite dimensional weight $U_q\mathfrak{g}$-module
$L(\lambda)$ with the highest weight $\lambda$:
$$
E_jv(\lambda)=0,\qquad K_j^\pm
v(\lambda)=q_j^{\pm\lambda_j}v(\lambda),\qquad
F_j^{\lambda_j+1}v(\lambda)=0.
$$
Such a weight $\lambda$ is said to be spherical if $L(\lambda)$ contains a
non-zero $U_q\mathfrak{k}$-invariant vector. This vector is unique up to a
constant multiple. The set $\Lambda$ of spherical weights is of the form
$\Lambda=\bigoplus\limits_{i=1}^r\mathbb{Z}_+\mu_i$, with
$\mu_1,\mu_2,\ldots,\mu_r$ being the so called fundamental spherical
weights. Here $r$ is the rank of the bounded symmetric domain $\mathbb{D}$.
In the simplest case $l=1$ the fundamental spherical weight $\mu$ is
$2\overline{\omega}$.

Equip the $U_q\mathfrak{g}$-module
$\mathbb{C}[X^\mathrm{spher}]_q\overset{\mathrm{def}}{=}
\bigoplus\limits_{\lambda\in\Lambda}L(\lambda)$ with a structure of
$U_q\mathfrak{g}$-module algebra in a way which is normally used in
producing quantum flag varieties. For any $\lambda',\lambda''\in\Lambda$,
$L(\lambda'+\lambda'')$ occurs in $L(\lambda')\otimes L(\lambda'')$ with
multiplicity 1. This allows one to introduce the morphisms of
$U_q\mathfrak{g}$-modules
$$
m_{\lambda',\lambda''}:L(\lambda')\otimes L(\lambda'')\to
L(\lambda'+\lambda''),\qquad m_{\lambda',\lambda''}:v(\lambda')\otimes
v(\lambda'')\mapsto v(\lambda'+\lambda''),
$$
and a structure of $U_q\mathfrak{g}$-module algebra in
$\mathbb{C}[X^\mathrm{spher}]_q$
$$
f'\cdot f''\overset{\mathrm{def}}{=}m_{\lambda',\lambda''}(f'\otimes
f''),\qquad f'\in L(\lambda'),\;f''\in L(\lambda'').
$$

Apply the Peter-Weyl decomposition
$\mathbb{C}[G]_q=\oplus_{\lambda\in\mathbb{Z}_+^l}(L(\lambda)\otimes
L(\lambda)^*)$ to obtain an embedding of $U_q\mathfrak{g}$-module algebras
$$
i:\mathbb{C}[X^\mathrm{spher}]_q\hookrightarrow\mathbb{C}[G]_q,\qquad
i:v(\lambda)\mapsto c_{\lambda,\lambda}^\lambda,\quad\lambda\in\Lambda.
$$
Here $c_{\lambda,\lambda}^\lambda$ are the matrix elements of
representation $\pi_\lambda$ associated to the $U_q\mathfrak{g}$-modules
$L(\lambda)$. This embedding allows one to treat the elements of
$\mathbb{C}[X^\mathrm{spher}]_q$ as q-analogues for sections of bundles on
the flag variety. Of course, $\mathbb{C}[X^\mathrm{spher}]_q$ has no zero
divisors, as these are absent in $\mathbb{C}[G]_q$.

Choose non-zero $U_q\mathfrak{k}$-invariants
$$\psi_j\in L(\mu_j),\qquad j=1,2,\ldots,r.$$
Of course, $\psi_1,\psi_2,\ldots,\psi_r\in\mathbb{C}[X^\mathrm{spher}]_q$.

\begin{proposition}
The elements $\psi_1,\psi_2,\ldots,\psi_r$ pairwise commute and the
multiplicative subset
$$
\Psi=\left\{\left.\psi_1^{j_1}\psi_2^{j_2}\cdots\psi_r^{j_r}\right|\:
j_1,j_2,\ldots,j_r\in\mathbb{Z}_+\right\}\subset
\mathbb{C}[X^\mathrm{spher}]_q
$$
is both right and left Ore set.
\end{proposition}

Let $\mathbb{C}[X^\mathrm{spher}]_{q,\Psi}$ be the localisation of
$\mathbb{C}[X^\mathrm{spher}]_q$ with respect to the multiplicative set
$\Psi$.

\begin{proposition}
There exists a unique extension of the structure of
$U_q\mathfrak{g}$-module algebra from $\mathbb{C}[X^\mathrm{spher}]_q$ onto
$\mathbb{C}[X^\mathrm{spher}]_{q,\Psi}$.
\end{proposition}

Our proof uses more elementary observations than those of a work by V.
Luntz and A. Rosenberg \cite{LR}.

Equip the algebra $\mathbb{C}[X^\mathrm{spher}]_{q,\Psi}$ with a
$\mathbb{Z}^r$-gradation so that
$$
\deg f=(j_1,j_2,\ldots,j_r)\;\Leftrightarrow\;f\in
L(j_1\mu_1+j_2\mu_2+\ldots+j_r\mu_r).
$$
The subalgebra of zero degree elements
$$
\mathbb{C}[\Omega]_q\overset{\mathrm{def}}{=}
\{f\in\mathbb{C}[X^\mathrm{spher}]_{q,\Psi}|\:\deg f=0\}
$$
is a $U_q\mathfrak{g}$-module algebra. This is just the q-analogue for the
algebra of regular functions on the open $K$-orbit $\Omega$.

In the simplest case $l=1$
$$i:\mathbb{C}[X^\mathrm{spher}]_q\hookrightarrow\mathbb{C}[SL_2]_q,$$
and the subalgebra $i(\mathbb{C}[X^\mathrm{spher}]_q)$ is generated by
\begin{equation}\label{gener}
t_{11}^2,\qquad t_{11}t_{12},\qquad t_{12}^2.
\end{equation}
In this case $r=1$, $i\psi_1=\mathrm{const}\cdot t_{11}t_{12}$, and $\Psi$
is an Ore set because $t_{11}t_{12}$ quasi-commutes with each of the
generators \eqref{gener}. In this special case $\mathbb{C}[\Omega]_q$ is
isomorphic to the algebra of Laurent polynomials in the indeterminate
$z=t_{12}^{-1}t_{11}$ just as in the classical case $q=1$. The action of
the generators $K^{\pm 1}$, $F$, $E$ of $U_q\mathfrak{sl}_2$ on Laurent
polynomials $f(z)$ is given by \eqref{KF_act}, \eqref{E_act}.

\bigskip

\stepcounter{section}

{\bf\thesection.} Many results related to the subjects of this talk remain
intact. For example, the works by H. Jakobsen and T. Tanisaki's team
\cite{J, T} study the $U_q\mathfrak{k}$-module algebras isomorphic to
$\mathbb{C}[\mathfrak{p}^-]_q$. The latter of those works suggests
q-analogues of Sato-Bernstein polynomials for quantum prehomogeneous vector
spaces in question. Note that such quantum prehomogeneous vector spaces
were found independently in \cite{T, J, SV}. q-Analogues for the Shilov
boundaries of the bounded symmetric domains can be produced. They can be
used to obtain geometric realizations for unitary principal degenerate
series of $U_q\mathfrak{g}$-modules (see \cite{Ber} for application of this
geometric realizations).

\end{document}